# A REPRESENTATION OF GIBBS MEASURE FOR THE RANDOM ENERGY MODEL

By Marie F. Kratz and Pierre Picco

*Université Paris V and CPT-FRUMAN*

In this work we consider a problem related to the equilibrium statistical mechanics of spin glasses, namely the study of the Gibbs measure of the random energy model. For solving this problem, new results of independent interest on sums of spacings for i.i.d. Gaussian random variables are presented. Then we give a precise description of the support of the Gibbs measure below the critical temperature.

**1. Introduction.** The problem of spin glasses is considered a challenge for the 21st century from both the theoretical physics point of view [47] and the mathematical point of view [46]. Spin glasses are alloys like Au–Fe that have a very small density of Fe. At very low temperature they present remarkable magnetic properties. On the experimental level, measurements of the magnetic susceptibility of such alloys at very low temperature were done by Cannela and Mydosh [8].

As mentioned in [8], the first model used to describe such an alloy was the Rudderman, Kittel, Kasuya and Yosida (RKKY) model, a long-range model with alternating sign in the interaction between the Fe atoms. Keeping just the alternating sign of the interaction, Edwards and Anderson introduced a model with short-range random coupling [21]. Despite a lack of rigorous results for the Edwards–Anderson model, a lot of theoretical progress has been made (see [36, 37, 48] and references therein). As a caricature of models with a spatial structure like the Ising model, mean field models such as the Curie–Weiss model have been introduced to give a simple partial description of physical phenomena, namely the existence of spontaneous magnetization at low temperature. The mean field model associated to the Edwards–Anderson model is called the Sherrington–Kirkpatrick (S–K) model [45]. From the theoretical physics point of view, much work had been

---

completed on the S–K model before Parisi [39, 40] introduced his famous replica symmetry breaking argument. A very important work by Mezard, Parisi and Virasoro [34] followed (see also [35]). Rigorous results on the S–K model can be found in [[1, 2, 10]–[12, 25, 28, 31, 46]]; however, we are still far from having complete results and this subject remains a very active field of research.

Due to the difficulty of studying the S–K model, even from a nonrigorous point of view, various mean field models have been introduced. Among them is the random energy model (REM) [16, 17], which is considered to be the simplest model of spin glasses (see [27]), and a whole class of models called generalized random energy models (GREM) that present a tree structure [[18]–[20]]. Basic rigorous results on the existence of the free energy for the REM and its fluctuations are given in [7, 22, 26, 38], whereas for the GREM, see [[9], [26]] and, more recently, [[5, 6]]. A very important fact is that the free energy of the GREM can also be found by using the replica symmetry breaking argument of Parisi (see [20]), since we know rigorously the explicit expression of the free energy (see [9]); this is an important test for the validity of the Parisi theory. Various pedagogical expositions on the REM and the GREM are available (e.g., [3, 4, 42, 44]).

The Hamiltonian or energy function of the REM can be written as

$$(1.1) \qquad H_{\mathrm{REM}}(\sigma) = -\frac{N^{1/2}}{2^{N/2}} \sum_{\alpha \subset \{1, \ldots, N\}} J_\alpha \sigma_\alpha,$$

where the sum is over all the $2^N$ subsets $\alpha$ of $\{1, \ldots, N\}$, $(J_\alpha, \alpha \subset \{1, \ldots, N\})$ is a family of i.i.d. standard Gaussians defined on a common probability space $(\Omega, \Sigma, \mathbb{P})$ and $\sigma_\alpha = \prod_{i \in \alpha} \sigma_i$ with $\sigma_\varnothing = 1$. The above way of writing the Hamiltonian of the REM emphasizes the mean field aspect of the model, since all spins interact with a strength that does not depend on their distances.

It is easy to see that the Hamiltonians $(H_{\mathrm{REM}}(\sigma), \sigma \in \{-1, +1\}^N)$ of the REM form a family of $2^N$ i.i.d. Gaussian random variables with mean 0 and variance $N$, defined on $(\Omega, \Sigma, \mathbb{P})$; see [24]. In particular, we can also write $(H_{\mathrm{REM}}(\sigma), \sigma \in \{-1, +1\}^N) = (-\sqrt{N} X_\sigma, \sigma \in \{-1, +1\}^N)$, where $(X_\sigma, \sigma \in \{-1, +1\}^N)$ is a family of i.i.d. Gaussian random variables with mean 0 and variance 1, or even $(-\sqrt{N} X_i, 1 \leq i \leq 2^N)$ when the relationship with the spins $\sigma$ is irrelevant. Here $\sqrt{N}$ ensures that the energy in a volume $N$ is $\mathbb{P}$-almost surely of order $N$, as required by statistical mechanics theory.

Given $\beta \geq 0$, the inverse temperature, let us denote by

$$(1.2) \qquad Z_N \equiv Z_N(\beta) = \sum_{\sigma \in \{-1, +1\}^N} e^{-\beta H_{\mathrm{REM}}(\sigma)}$$

the finite volume partition function and by

$$(1.3) \qquad F_N(\beta) = -\frac{1}{\beta N} \log Z_N(\beta)$$

the finite volume free energy.

It was proved in [38] that for all $\beta \geq 0$, $\lim_{N \to \infty} F_N(\beta) = F(\beta)$ exists $\mathbb{P}$-almost surely and in $L^p(\Omega, \mathbb{P})$ for $1 \leq p < \infty$ (note that what is called free energy in [38] we denote $-\beta F_N(\beta)$). The nonrandom function $F(\beta)$ is twice differentiable in $\beta$, whose second derivative has a jump at $\beta_c := \sqrt{2 \log 2}$. In fact, $F(\beta)$ is equal to $-\beta/2 - \beta_c^2/(2\beta)$ for $\beta < \beta_c$ and $-\beta_c$ for $\beta \geq \beta_c$, as expected from the results of [16]. In the physics literature this is called a third order phase transition, as mentioned explicitly in the title of the first rigorous paper on the REM (see [22]). More results on the asymptotics of $F_N(\beta)$ can be found in [7, 26, 38].

What is important to note is that when $\beta > \beta_c$, the main contribution to $F_N(\beta)$ comes from the terms in $Z_N(\beta)$ that have the lowest possible energies, which are the $\sigma$ that have $H_{\mathrm{REM}}(\sigma)$ of order $-N\sqrt{2 \log 2} \equiv -N\beta_c$, that is, the ($\mathbb{P}$-almost sure) asymptotic value of the minimum of $2^N$ Gaussian random variables with mean 0 and variance $N$; this fact is explicitly mentioned in all pedagogical expositions (e.g., [3, 4, 42]). In the physics literature, such a system is said to be frozen in the sense that, in the whole range of $\beta \geq \beta_c$, the lowest possible energies give the main contribution to the free energy, while in a "nonfrozen" system, this contribution occurs only at zero temperature. This can be seen on $F(\beta)$ since when $\beta > \beta_c$, the free energy does not depend on $\beta$.

To study the fluctuations of the free energy $F_N(\beta)$, we consider the number of random variables of the sample that are below some well chosen nonrandom energy, and then use classical convergence to a Poisson point process [30] to describe some fluctuations of the model. This was done in [26] for the Boltzmann factor $\exp(-\beta H_{\mathrm{REM}}(\sigma))$ (see also [3]) and in [7] for the partition function.

The aim of this article is to study the finite volume Gibbs measure $\mu_{N,\beta}$ of the REM when $\beta > \beta_c$ and $N$ is large enough. Here $\mu_{N,\beta}$ is defined for each $N$ as the random probability measure on $\{-1, +1\}^N$ which gives the configuration $\sigma$ the weight

$$(1.4) \qquad \mu_{N,\beta}(\sigma) \equiv \frac{e^{-\beta H_{\mathrm{REM}}(\sigma)}}{Z_N},$$

where $Z_N$ is defined as in (1.2).

For a given sample $(H_{\mathrm{REM}}(\sigma), \sigma \in \{-1, +1\}^N)$, when $\beta > \beta_c$, the main question is: What are the sample: dependent configurations $\sigma$, where the Gibbs measure is concentrated, when $N$ is very large?

Note first that the sample-dependent configuration $\sigma^{(1)}$ which corresponds to the minimal value of the sample $H_{\text{REM}}(\sigma)$ is clearly among these configurations. So instead of considering the variables that are below a nonrandom energy as was done to get the Poisson point process mentioned above, it is better to consider the variables that are above the minimal one that is a sample-dependent energy. It is clear that the order statistics of the random variables $H_{\text{REM}}(\sigma)$, namely

$$(1.5) \quad H_{\text{REM}}(\sigma^{(2^N)}) \geq H_{\text{REM}}(\sigma^{(2^{N-1})}) \geq \cdots \geq H_{\text{REM}}(\sigma^{(2)}) \geq H_{\text{REM}}(\sigma^{(1)}),$$

come into play. To obtain information on the support of the Gibbs measure, we can subtract the minimal energy from all the Hamiltonians; therefore, the spacings $H_{\text{REM}}(\sigma^{(k+1)}) - H_{\text{REM}}(\sigma^{(k)})$ and the sums of spacings $H_{\text{REM}}(\sigma^{(k)}) - H_{\text{REM}}(\sigma^{(1)})$ are the basic objects we need to study.

Since our main objective is more the study of the Gibbs measure of the REM than the sum of spacings per se, we first present a list of questions that we asked ourselves when we started this work. Then the strategy we used and the results stated in the next section about the sums of spacings will be clearer.

First, how many successive terms $k = k_N$ in the sum of spacings do we need to take to have a "good" approximation of the Gibbs measure when considering the probability measure on $\{-1, +1\}^N$ that havs only these $k$ terms?

There are two points to be answered precisely: What does "good approximation" mean and in what $\mathbb{P}$-probability sense can we expect that such approximation holds? For the first point we choose the total variation distance for probability measures on $\{-1, +1\}^N$. Concerning the $\mathbb{P}$-probabilistic sense, some care is needed. We first note that there is a scaling factor $\sqrt{N}$ in the definition of the Hamiltonian, as mentioned previously. Recalling that we have a sample of $2^N$ Gaussian random variables, this factor can be written, up to a constant, as $\sqrt{2 \log 2^N}$. So the first question is related to the behavior of $\sqrt{2 \log 2^N}[H_{\text{REM}}(\sigma^{(k)}) - H_{\text{REM}}(\sigma^{(1)})]$ as a function of $k$ for $k \equiv k_N$ that diverges with $N$. Since we have chosen the total variation distance, it will be sufficient to neglect all the terms in the Gibbs measures that are larger than $k_N$. An easy way to do this is to require that $\sum_{k=k_N}^{2^N} \exp -\{(\beta/\beta_c)(\sqrt{2 \log 2^N}[H_{\text{REM}}(\sigma^{(k)}) - H_{\text{REM}}(\sigma^{(1)})])\}$ goes to zero. The strongest $\mathbb{P}$-probability sense we can expect to get is $\mathbb{P}$-almost surely. At this point, it is important to note that we do not expect that the Gibbs measure itself, or its approximation converges almost surely since they merely converge in law. However, their total variation distance, being a difference, could perfectly converge almost surely to 0. A similar fact was proved not for the Gibbs measure as a measure on $\{-1, +1\}^N$, but for an induced measure in the random field Curie–Weiss model [**?** ]. The worst that could happen by taking such a strong convergence would be that we get a $k_N$ that is rather large.

The main result we need for the sum of spacings is given in Proposition 2.2. It is stated for a sample of $n$ Gaussian random variables, instead of $2^N$, since these results have independent interest.

The second question is: Is the approximated Gibbs measure, that is, the one with only $k_N$ terms obtained previously, the uniform measure on $k_N$ points chosen without replacement within $2^N$ points, a point mass at the minimum or some other value? To be able to answer such questions, the relative weights of the approximated Gibbs measure and, therefore, the first $k_N$ sum of spacings come into play. We use a classical representation of the Gaussian random variables in term of uniform random variables, mainly because a lot of explicit distributions for the spacings of uniform random variables are available. The results needed for the sum of spacings are stated in Proposition 2.3, where the first $k_N$ sums of spacings can be represented as successive partial weighted sums of exponential random variables.

Note also that since the $k_N$ was found to be able to neglect some tails of the Gibbs measure, $\mathbb{P}$-almost surely, it could be that the Gibbs measure has most of its mass concentrated on the $k_0$ spin configuration with $k_0$ that does not depend on $N$. We will show that this is not the case in general. Results for the approximated Gibbs measure and the Gibbs measure are given in Section 3.

Another question is related to the $\beta$-dependence of the Gibbs measure. As mentioned above, the system is frozen at the level of the free energy. Does a similar fact hold for the Gibbs measure? When the limit $\beta \uparrow \infty$ is taken after the limit $N \uparrow \infty$, the Gibbs measure converges to a point mass at the spin configuration that realizes the absolute minimum as is proved in Section 3. So to exclude the possibility that the system is frozen at the level of the Gibbs measure, we estimate from below the total variation distance between the Gibbs measure at finite $\beta > \beta_c$ and its limit when $\beta \uparrow \infty$.

Assuming that some or all of the foregoing questions have been resolved—in particular, that the number $k(N)$ has been determined to obtain, $\mathbb{P}$-almost surely a good approximation of the Gibbs measure—then is there an "easy" way to construct this approximated Gibbs measure? This is the subject of the last section.

## 2. Some results involving order statistics and sums of spacings.

2.1. *Notation and recollections.* In this section, $(X_i, i \in \{1, \ldots, n\})$ are i.i.d. standard Gaussian random variables with distribution function $\Phi$, $(U_i, i \in \{1, \ldots, n\})$ are i.i.d. $\mathcal{U}(0, 1)$, random variables and $(W_i, i \in \{1, \ldots, n\})$ are i.i.d. exponentially distributed r.v. $\mathcal{E}(1)$. We define $S_0 := 0$ and $S_m := \sum_{k=1}^{m} W_k$, $m \in \mathbb{N}^*$. In general we denote by $Y_{1,n} \leq \cdots \leq Y_{n,n}$ the order statistics associated to some random variables $(Y_i, i \in \{1, \ldots, n\})$. We set $U_{0,n} := 0$ and $U_{n+1,n} := 1$.

Let $G$ be defined by $1 - \Phi(G(u)) = u$, $0 < u < 1$. Then $G$ satisfies (see [14], page 818), as $u \downarrow 0$,

$$(2.1) \quad G(u) = \sqrt{2 \log\left(\frac{1}{u}\right)} - \frac{\log \log(1/u) + \log(4\pi)}{2\sqrt{2 \log(1/u)}} + O\left(\frac{(\log \log(1/u))^2}{(\log(1/u))^{3/2}}\right).$$

As it is standard, we can construct the Gaussian random variables by using $X_i = G(U_i)$, and since $G$ is decreasing, we have $X_{i,n} = G(U_{n-i+1,n})$. Then, by symmetry, we have the identity in distribution (denoted $\overset{d}{=}$)

$$(2.2) \quad X_{i,n} - X_{1,n} \overset{d}{=} G(U_{1,n}) - G(U_{i,n}).$$

Recall the results

$$(2.3) \quad \{U_{i,n}, 0 \leq i \leq n+1\} \overset{d}{=} \{S_i/S_{n+1}, 0 \leq i \leq n+1\}$$

and

$$(2.4) \left\{\xi_i^{(n)} := \left(\frac{U_{i,n}}{U_{i+1,n}}\right)^i; 1 \leq i \leq n\right\} \quad \text{are i.i.d. } \mathcal{U}(0,1) \text{ random variables,}$$

which can be found in [15], Lemmas 1 and 2, or in [43] and [32], respectively.

2.2. *New results.* Before stating and proving the main result of this section, Proposition 2.2, we present some preliminary results:

LEMMA 2.1. *For all $\delta > 0$, there exists $n_0 = n_0(\delta)$ such that for all $n > n_0$, we have*

$$(2.5) \quad \mathbb{P}\left[\left|\log\left(\frac{1}{U_{1,n}}\right) - \log n\right| \leq 2\log(\log n)^{1+\delta}\right] \geq 1 - \frac{4}{(\log n)^{1+\delta}}.$$

PROOF. Equation (2.3) implies that

$$
\begin{aligned}
&\mathbb{P}\left[\left|\log\left(\frac{1}{U_{1,n}}\right) - \log n\right| \leq 2\log(\log n)^{1+\delta}\right] \\
(2.6) \quad &\geq \mathbb{P}[|\log S_1| \leq \log(\log n)^{1+\delta}] \\
&\quad - \mathbb{P}\left[\left|\log\left(\frac{S_{n+1}}{n+1}\right) + \log\left(1 + \frac{1}{n}\right)\right| \geq \log(\log n)^{1+\delta}\right].
\end{aligned}
$$

Since $S_1$ is exponentially distributed, we obtain

$$(2.7) \quad \mathbb{P}[|\log S_1| \leq \log(\log n)^{1+\delta}] \geq 1 - \exp\{-(\log n)^{1+\delta}\} - \frac{1}{(\log n)^{1+\delta}}.$$

By using the exponential Markov inequality and the fact that $(1-x)^{-1} \leq e^{x+x^2}$ for $|x| \leq 1/2$, we get that $\forall\, 0 < \varepsilon < 1$, $\forall\, n \geq 1$,

$$(2.8) \qquad \mathbb{P}\left[\left|\frac{S_n}{n} - 1\right| \geq \varepsilon\right] \leq 2e^{-n\varepsilon^2/4},$$

which implies after a short computation that, $\forall\, 0 < \varepsilon < 1$, $\forall\, \delta > 0$, $\exists\, n_{\varepsilon,\delta} = \exp(2/(1-\varepsilon))^{1/(1+\delta)}$ such that $\forall\, n > n_{\varepsilon,\delta}$,

$$(2.9) \qquad \mathbb{P}\left[\left|\log\left(\frac{S_{n+1}}{n+1}\right) + \log\left(1 + \frac{1}{n}\right)\right| \geq \log(\log n)^{1+\delta}\right] \leq 2e^{-n\varepsilon^2/4}.$$

Then combining (2.6), (2.7) and (2.9) and taking $\varepsilon = 1/2$ entails the existence of $n_0 \equiv n_0(\delta)$ such that (2.5) holds. $\quad\square$

PROPOSITION 2.1. *We have*

$$(2.10) \qquad \begin{aligned} \log\left(\frac{U_{1,n}}{U_{j,n}}\right) &\stackrel{d}{=} -\sum_{i=1}^{j-1} \frac{W_i}{i} \\ &\stackrel{d}{=} -Z + \sum_{i=j}^{\infty}\left(\frac{W_i}{i} - \frac{E[W_i]}{i}\right) - \sum_{i=1}^{j-1} \frac{E[W_i]}{i}, \end{aligned}$$

*where $Z$ is a random variable such that $\forall\, x > 0$,*

$$(2.11) \quad P[Z \geq x] \leq e^{-(x\sqrt{3}/2)+15/4} \quad and \quad P[Z \leq -x] \leq e^{-(x\sqrt{3}/2)+15/4}.$$

*Moreover, for all $j > e$, we have*

$$(2.12) \quad \mathbb{P}\left[\log\left(\frac{U_{1,n}}{U_{j,n}}\right) \leq -Z + \sqrt{\frac{3\log j}{j}} - \log j\right] \geq 1 - \frac{\exp(3((\log j)^2/j))}{j^{3/4}}.$$

*More generally, $\forall\, c > 0$, $\forall\, 0 < \varepsilon < 1$, $\exists\, j_0 = j_0(c, \varepsilon)$ and $\forall\, j \geq j_0$,*

$$(2.13) \qquad \mathbb{P}\left[\log\left(\frac{U_{1,n}}{U_{j,n}}\right) \leq -Z + \sqrt{\frac{c\log j}{j}} - \log j\right] \geq 1 - \frac{1}{j^{(1-\varepsilon)c/4}}.$$

The proof of this proposition necessitates the following result:

LEMMA 2.2. *For all positive $x$, for all positive integers $j_0$, $j_1$ (with $j_1 \geq j_0$) and for all positive $t$ such that $t^2/j_0^2 \leq 3/4$, we have*

$$(2.14) \qquad \begin{aligned} &\mathbb{P}\left[\sum_{l=j_0}^{j_1} \frac{W_l - E[W_l]}{l} \geq x\right] \\ &\qquad \leq \exp\left\{-tx + \frac{t^2}{j_0}\left(1 + 2\frac{t^2}{j_0^2}\right)\left(1 + \frac{1}{j_0}\right)\right\}. \end{aligned}$$

*In particular, for $t^2 \leq 3/4$,*

$$(2.15) \qquad \mathbb{P}\left[\sum_{l=1}^{j} \frac{W_l - E[W_l]}{l} \geq x\right] \leq e^{-tx + 2t^2(1+2t^2)}.$$

PROOF. The proof is a consequence of exponential Markov inequality and Jensen inequality. Indeed, we can write

$$\mathbb{P}\left[\sum_{l=j_0}^{j_1} \frac{W_l - E[W_l]}{l} \geq x\right]$$

$$\leq e^{-tx} \prod_{l=j_0}^{j_1} E\left(\exp\left\{t\frac{W_l - E[W_l']}{l}\right\}\right)$$

$$\leq e^{-tx} \prod_{l=j_0}^{j_1} E\left(\exp\left\{t\frac{W_l - W_l'}{l}\right\}\right)$$

$$= e^{-tx} \prod_{l=j_0}^{j_1} \frac{1}{1 - t^2/l^2}$$

$$\leq \exp\left\{-tx + \left(1 + \frac{2t^2}{j_0^2}\right)t^2 \sum_{l=j_0}^{j_1} \frac{1}{l^2}\right\},$$

with $(W_l')$ i.i.d. standard exponentials, that are independent of $(W_l)$ and by using the inequality $(1-x)^{-1} \leq e^{x+2x^2}$ for $0 \leq x \leq 3/4$ at the last step. Now (2.14) follows since

$$\sum_{l=j_0}^{j_1} \frac{1}{l^2} \leq \frac{1}{j_0^2} + \int_{j_0}^{j_1} \frac{dx}{x^2} \leq \frac{1}{j_0}\left(1 + \frac{1}{j_0}\right). \qquad \square$$

PROOF OF PROPOSITION 2.1. We have, by using (2.4),

$$(2.16) \qquad \log\left(\frac{U_{1,n}}{U_{j,n}}\right) = \log\left(\frac{U_{1,n}}{U_{2,n}} \cdots \frac{U_{j-1,n}}{U_{j,n}}\right)$$

$$= \sum_{i=1}^{j-1} \frac{1}{i} \log\left(\frac{U_{i,n}}{U_{i+1,n}}\right)^i = \sum_{i=1}^{j-1} \frac{1}{i} \log \xi_i^{(n)},$$

from which we deduce the first equality of (2.10). Now we can write

$$\sum_{i=1}^{j-1} \frac{W_i}{i} = \sum_{i=1}^{\infty} \frac{W_i - E[W_i]}{i} - \sum_{i=j}^{\infty} \frac{W_i - E[W_i]}{i} + \sum_{i=1}^{j-1} \frac{E[W_i]}{i},$$

where the two series converge not only in quadratic mean, but also almost surely by the Lévy theorem.

Calling $Z = \sum_{i=1}^{\infty}(W_i - E[W_i])/i$, the tail of the distribution of $Z$ decreases exponentially, as mentioned in (2.11). Indeed, it suffices, on one hand, to apply (2.15) for $t = \sqrt{3}/2$ to get the first property on $Z$ and, on the other hand, to note that since

$$\mathbb{P}\left[\sum_{l=1}^{\infty}\frac{W_l - E[W_l]}{l} \leq -x\right] = \mathbb{P}\left[\sum_{l=1}^{\infty}\frac{-W_l + E[W_l]}{l} \geq x\right],$$

by making exactly the same computation as in Lemma 2.2, we get the same estimate and therefore the second property on $Z$.

By applying (2.14) with $x = \sqrt{(c\log j)/j}$ and $t = \sqrt{cj\log j}/2$, we can say that, for all $c > 0$, for all $j > 1$ such that $3j/\log j \geq c$,

$$\mathbb{P}\left[\sum_{l=j}^{\infty}\frac{W_l - E[W_l]}{l} \geq \sqrt{\frac{c\log j}{j}}\right] \leq \exp\left\{-\frac{c}{4}\log j\left(1 - \frac{1 + c\log j}{j}\right)\right\}.$$

Therefore, for all $c > 0$, for all $0 < \varepsilon < 1$ and $\exists j_0 = j_0(c, \varepsilon)$ such that $j \geq j_0$, we have $(1 + c\log j)/j \leq \varepsilon$, which entails (2.13), whereas for $j > e$, taking $c = 3$ gives (2.12). $\square$

Now we can state the main result of this section.

PROPOSITION 2.2. *For all $0 < \delta$, for all $0 < \varepsilon < 1$ and for all $k_n$ satisfying*

$$(2.17) \qquad k_n \uparrow \infty, \qquad \frac{k_n}{\log n} \uparrow \infty \quad and \quad \frac{\log k_n}{\log n} \downarrow 0 \qquad as \ n \to \infty,$$

*by defining, for $0 < \lambda < 1$ and $0 < \alpha < 1$,*

$$(2.18) \qquad \lambda_n := \lambda\frac{\log n}{n^{\alpha}}$$

*and*

$$(2.19) \qquad \tilde{\lambda}_n := \lambda_n\left(1 + 2\sqrt{\frac{2}{\lambda}}\frac{(1 - \lambda_n)}{n^{(1-\alpha)/2}}\right),$$

*there exists $n_0 = n_0(\varepsilon, \delta, \lambda, \alpha)$ such that $\forall\, n \geq n_0$, there exists $\Omega_n \subset \Omega$, with*

$$(2.20) \qquad \mathbb{P}[\Omega_n] \geq 1 - 4\left(\frac{1}{(\log n)^{1+\delta}} + \frac{e^{-k_n/16}}{1 - e^{-1/16}}\right),$$

*such that on $\Omega_n$, for all $j$ such that $k_n \leq j \leq n\tilde{\lambda}_n$, we have*

$$(2.21) \qquad \sqrt{2\log n}(X_{j,n} - X_{1,n}) \geq 2\log n\left(1 - \sqrt{1 - \frac{\log j}{\log n}}\right)(1 - \varepsilon),$$

*while for all $j$ such that $n\tilde{\lambda}_n \leq j \leq n$,*

$$(2.22) \qquad \sqrt{2\log n}(X_{j,n} - X_{1,n}) > 2\log n\sqrt{1-\varepsilon} - G(\lambda_n)\sqrt{2\log n}.$$

*In particular, on $\Omega_n$, for $\log j/\log n \uparrow 1$,*

$$(2.23) \qquad \sqrt{2\log n}(X_{j,n} - X_{1,n}) > 2\log n\sqrt{1-\varepsilon} - G(\lambda_n)\sqrt{2\log n}.$$

REMARKS.

1. If $\lambda_n$ is chosen as a constant $\lambda$ independent of $n$, $0 < \lambda < 1$, (2.22) gives that $\forall j > n\lambda$, $\sqrt{2\log n}(X_{j,n} - X_{1,n}) \geq 2\log n\sqrt{1-\varepsilon} - G(\lambda)(\sqrt{2\log n}) > 2\log n(1-\varepsilon)$ if $n > n_0(\lambda,\varepsilon)$ for some $n_0(\lambda,\varepsilon)$. When entering various regimes as $n\tilde{\lambda}_n^{(1)} < j < n\tilde{\lambda}_n^{(2)}$ with $\lambda_n^{(i)} \downarrow 0$ for $i = 1, 2$, a cancellation could occur between the two terms in the right-hand side of (2.22). The choice of $\lambda_n$ in (2.18) then allows us to see such cancellation, since for $\alpha < 1/4$ it provides that, $\forall j \geq n\tilde{\lambda}_n$, on $\Omega_n$,

$$(2.24) \quad \sqrt{2\log n}(X_{j,n} - X_{1,n}) \geq 2(1 - \sqrt{\alpha})(1 - \varepsilon)\log n \geq (1 - \varepsilon)\log n.$$

2. Note that the lower bounds in (2.21) and (2.22), even when obtained by two completely different methods, are of the same order $2(1 - \sqrt{\alpha})\log n$ when $j = n\tilde{\lambda}_n$.

3. Note that for $n = 2^N$, under (2.17), we have $\sum_{N=1}^{\infty} \mathbb{P}[\Omega_{2^N}^c] < \infty$. In the next section, this allows us to get some results that are true $\mathbb{P}$-almost surely for all but a finite number of indices $N$ by the first Borel–Cantelli lemma.

4. For completeness, even though it will not be used in the next section, we show that for $0 < \eta < 1$, there exists $\Omega_{k_n}^*$ with $\mathbb{P}(\Omega_{k_n}^*) \geq 1 - 2\exp\{-\sqrt{3}/2(\log k_n)^\eta\}$, such that on $\Omega_{k_n}^*$, for all $j$ such that $k_n \leq j \leq n\tilde{\lambda}_n$ and $\log j/\log n \downarrow 0$, (2.21) can be refined as

$$(2.25) \qquad \sqrt{2\log n}(X_{j,n} - X_{1,n}) \geq \log j\left(1 - \frac{1}{(\log j)^{1-\eta}}\right).$$

PROOF OF PROPOSITION 2.2. Let us separately consider the two cases as in (2.21) and (2.22).

CASE (i). $j$ is such that $n\tilde{\lambda}_n \leq j \leq n$. Let us denote by

$$(2.26) \quad \mathcal{L} = \mathcal{L}(n, \tilde{\lambda}_n, \varepsilon) := \bigcup_{j=n\tilde{\lambda}_n}^{n} \{\sqrt{2\log n}(X_{n,n} - X_{n-j+1,n}) < f(n,\lambda)\}$$

with $f(n,\lambda) := 2\log n\sqrt{1-\varepsilon} - \sqrt{2\log n}G(\lambda_n)$.

Note that in this case we are working where $j$ is large with $(X_{n,n} - X_{n-j+1,n})$ instead of $(X_{j,n} - X_{1,n})$ since these spacings are equal in distribution. We claim that

$$(2.27) \qquad \mathbb{P}[\mathcal{L}] \leq 2e^{-4(1-\lambda_n)\log n} + \mathbb{P}[X_{n,n} < \sqrt{2(1-\varepsilon)\log n}\,],$$

which entails that $\forall\, 0 < \lambda < 1$, $\exists\, n_0$ and $\forall\, n \geq n_0$,

$$(2.28) \qquad \mathbb{P}[\mathcal{L}] \leq \frac{2}{n^2} + \mathbb{P}[X_{n,n} < \sqrt{2(1-\varepsilon)\log n}\,],$$

after having noticed that $\forall\, \lambda$, $0 < \lambda < 1$, $\forall\, \alpha$, $0 < \alpha < 1$, $\exists\, n_0 = n_0(\lambda, \alpha)$, s.t. $\forall\, n \geq n_0$, $4(1 - \lambda_n) \geq 2$.

Now, (2.27) is an immediate consequence of

$$(2.29) \qquad \mathbb{P}[\mathcal{L} \cap \{X_{n,n} \geq \sqrt{2(1-\varepsilon)\log n}\,\}] \leq 2e^{-4(1-\lambda_n)\log n}.$$

To prove (2.29), note that, on one hand,

$$(2.30) \qquad \begin{aligned} & \mathcal{L} \cap \{X_{n,n} \geq \sqrt{2(1-\varepsilon)\log n}\,\} \\ & \quad \subset \widetilde{\mathcal{L}} := \bigcup_{j=n\widetilde{\lambda}_n}^{n} \{\sqrt{2\log n}(\sqrt{(1-\varepsilon)2\log n} - X_{n-j+1,n}) < f(n,\lambda)\} \end{aligned}$$

and, on the other hand,

$$(2.31) \qquad \mathbb{P}[\mathcal{L} \cap (X_{n,n} \geq \sqrt{2(1-\varepsilon)\log n}\,)] \leq \mathbb{P}[\mathcal{A}] + \mathbb{P}[\widetilde{\mathcal{L}} \cap \mathcal{A}^c]$$

with

$$(2.32) \qquad \begin{aligned} \mathcal{A} := \Bigg\{ & \left| \sum_{i=1}^{n} (\mathbb{1}_{\{X_i \leq G(\lambda_n)\}} - \mathbb{E}[\mathbb{1}_{\{X_i \leq G(\lambda_n)\}}]) \right| \\ & \geq \lambda_n(1-\lambda_n)n^{(1+\alpha)/2}2\sqrt{2/\lambda} \Bigg\}, \end{aligned}$$

where $\mathbb{E}[\mathbb{1}_{\{X_i \leq G(\lambda_n)\}}] = 1 - \lambda_n$ by definition of $G$.

To estimate $\mathbb{P}[\mathcal{A}]$, we use Bernstein's inequality [41], namely:

LEMMA 2.3. *Let* $(Y_i)_i$ *be independent random variables, mean* 0*, s.t.* $\mathbb{E}[Y_i^2] < \infty$ *and* $|Y_i| \leq 1, \forall\, i$. *If*

$$D_n := \sum_{i=1}^{n} \mathbb{E}[Y_i^2]$$

*and if* $0 < t < \sqrt{D_n}$, *then*

$$\mathbb{P}\left[ \left| \sum_{i=1}^{n} Y_i \right| \geq t\sqrt{D_n} \right] \leq 2e^{-t^2/2}.$$

By using this inequality for $Y_i := \mathbb{1}_{\{X_i \leq G(\lambda_n)\}} - \mathbb{E}[\mathbb{1}_{\{X_i \leq G(\lambda_n)\}}]$ with $t = 2n^{\alpha/2}\sqrt{\lambda_n(1-\lambda_n)2/\lambda}$, since $D_n = n\lambda_n(1-\lambda_n)$, $\exists n_0(\lambda, \alpha)$, $\forall n \geq n_0$, $t < \sqrt{D_n}$, therefore we get that $\forall n \geq n_0$,

$$(2.33) \qquad\qquad P[\mathcal{A}] \leq 2e^{-4(1-\lambda_n)\log n}.$$

On the other hand,

$$(2.34) \quad \mathcal{A}^c \subset \left\{ n(1-\lambda_n) - \lambda_n(1-\lambda_n)n^{(1+\alpha)/2}2\sqrt{2/\lambda} < \sum_{i=1}^n \mathbb{1}_{\{X_i \leq G(\lambda_n)\}} \right\},$$

which implies that on $\mathcal{A}^c$ the number of random variables $X_i$ less than $G(\lambda_n)$ is greater than $1 + n(1-\tilde{\lambda}_n)$; hence, $\mathcal{A}^c \subset \{X_{1+n(1-\tilde{\lambda}_n),n} \leq G(\lambda_n)\}$ and so

$$(2.35) \quad \mathbb{P}\left[ \left( \bigcup_{k=1}^{n-n\tilde{\lambda}_n+1} X_{k,n} > G(\lambda_n) \right) \cap (X_{1+n-n\tilde{\lambda}_n,n} \leq G(\lambda_n)) \right] = 0.$$

Combining (2.31), (2.33) and (2.35) entails (2.27).

Now (2.28) and the fact that for all $0 < \varepsilon < 1$,

$$(2.36) \quad \mathbb{P}[X_{n,n} < \sqrt{2(1-\varepsilon)\log n}] = \Phi^n(\sqrt{2(1-\varepsilon)\log n}) < e^{-n^\varepsilon/\sqrt{2(1-\varepsilon)\log n}}$$

give

$$(2.37) \qquad\qquad \mathbb{P}[\mathcal{L}] \leq \frac{2}{n^2} + e^{-n^\varepsilon/\sqrt{2(1-\varepsilon)\log n}},$$

which leads to (2.22), since, for all $0 < \varepsilon < 1$, for all $\delta > 0$, for all $0 < \lambda < 1$, for all $0 < \alpha < 1$, for all $k_n$ satisfying (2.17), $\exists n_0 = n_0(\varepsilon, \delta, \lambda, \alpha)$ such that $\forall n \geq n_0$,

$$(2.38) \qquad \frac{2}{n^2} + e^{-n^\varepsilon/\sqrt{2(1-\varepsilon)\log n}} \ll 2\left( \frac{1}{(\log n)^{1+\delta}} + \frac{e^{-k_n/16}}{1-e^{-1/16}} \right).$$

CASE (ii).  $k_n \leq j \leq n\tilde{\lambda}_n$. Note that, because of (2.1), we have, for $u \downarrow 0$, $v \downarrow 0$,

$$G(u) - G(v) = \sqrt{2}\left( \sqrt{\log\left(\frac{1}{u}\right)} - \sqrt{\log\left(\frac{1}{v}\right)} \right)$$
$$+ O\left( \frac{\log\log(1/u)}{\sqrt{\log(1/u)}} \right) + O\left( \frac{\log\log(1/v)}{\sqrt{\log(1/v)}} \right).$$

For $\gamma > 0$, let

$$\Omega_{\tilde{\lambda}_n} := \left\{ \sup_{1 \leq j \leq n\tilde{\lambda}_n} U_{j,n} = U_{n\tilde{\lambda}_n,n} \leq (1+\gamma)\tilde{\lambda}_n \right\}.$$

Since

$$\mathbb{P}[\Omega_{\tilde{\lambda}_n}^c] \leq \mathbb{P}\left[\Omega_{\tilde{\lambda}_n}^c \cap \left(\left|\frac{S_{n+1}}{n+1} - 1\right| < \varepsilon\right)\right] + \mathbb{P}\left[\left|\frac{S_{n+1}}{n+1} - 1\right| \geq \varepsilon\right]$$

$$\leq \mathbb{P}\left[\frac{S_{n\tilde{\lambda}_n}}{S_{n+1}} \geq (1+\gamma)\tilde{\lambda}_n, S_{n+1} \geq (n+1)(1-\varepsilon)\right] + \mathbb{P}\left[\left|\frac{S_{n+1}}{n+1} - 1\right| \geq \varepsilon\right]$$

$$\leq \mathbb{P}\left[\frac{S_{n\tilde{\lambda}_n}}{n\tilde{\lambda}_n} \geq (1+\gamma)(1-\varepsilon)\frac{n+1}{n}\right] + \mathbb{P}\left[\left|\frac{S_{n+1}}{n+1} - 1\right| \geq \varepsilon\right],$$

by using (2.3) in the second inequality, it follows from (2.8) with $\varepsilon = \gamma/(2+\gamma)$, that $\mathbb{P}[\Omega_{\tilde{\lambda}_n}^c] \leq 4\exp(-n\tilde{\lambda}_n(\gamma^2)/(4(2+\gamma)^2))$.

On $\Omega_{\tilde{\lambda}_n}$, we have

$$G(U_{1,n}) - G(U_{j,n})$$

$$(2.39) \qquad = \sqrt{2\log\left(\frac{1}{U_{1,n}}\right)}\left(1 - \sqrt{1 + \frac{\log(U_{1,n}/U_{j,n})}{\log(1/U_{1,n})}}\right)$$

$$+ O\left(\frac{\log\log(1/\tilde{\lambda}_n)}{\sqrt{\log(1/\tilde{\lambda}_n)}}\right).$$

Let $\Omega_{n,\delta} \subset \Omega$ be defined by

$$\Omega_{n,\delta} := \left\{\left|\frac{\log(1/U_{1,n})}{\log n} - 1\right| \leq 2\frac{\log(\log n)^{1+\delta}}{\log n}\right\}.$$

Then Lemma 2.1 gives

$$(2.40) \qquad P[\Omega_{n,\delta}] \geq 1 - \frac{4}{(\log n)^{1+\delta}}.$$

Let $\tilde{\Omega}_{n,\delta} := \Omega_{n,\delta} \cap \Omega_{\tilde{\lambda}_n}$. Then by combining (2.2), (2.5) and (2.39), on $\tilde{\Omega}_{n,\delta}$ we obtain that

$$\sqrt{2\log n}(X_{j,n} - X_{1,n})$$

$$(2.41) \qquad = 2\log n\left(1 - \sqrt{1 + \frac{\log(U_{1,n}/U_{j,n})}{\log n}}\right)\left(1 + O\left(\frac{\log\log n}{\sqrt{\log n}}\right)\right).$$

Again using (2.3), we can write

$$(2.42) \qquad \frac{\log(U_{1,n}/U_{j,n})}{\log n} = \frac{\log(S_1)}{\log n} - \frac{\log(S_j/j)}{\log n} - \frac{\log j}{\log n}.$$

Let

$$(2.43) \qquad \tilde{\Omega}_j := \left\{\left|\frac{S_j}{j} - 1\right| < \frac{1}{2}\right\} \quad \text{and} \quad \Omega_{k_n} := \left(\bigcap_{j=k_n}^{n\tilde{\lambda}_n} \tilde{\Omega}_j\right) \cap \Omega_{n,\delta}.$$

Then by using (2.40) and (2.8) with $\varepsilon = 1/2$, we get

$$(2.44) \qquad \mathbb{P}[\Omega_{k_n}^c] \leq 4\left(\frac{1}{(\log n)^{1+\delta}} + \frac{e^{-k_n/16}}{1 - e^{-1/16}}\right).$$

On $\Omega_{k_n}$, we have

$$(2.45) \qquad \frac{\log(U_{1,n}/U_{j,n})}{\log n} = -\frac{\log j}{\log n} + O\left(\frac{\log(\log n)^{1+\delta}}{\log n}\right),$$

which combined with (2.41) gives

$$(2.46) \qquad \begin{aligned} &\sqrt{2\log n}\,(X_{j,n} - X_{1,n}) \\ &= 2\log n\left(1 - \sqrt{1 - \frac{\log j}{\log n} + O\left(\frac{\log(\log n)^{1+\delta}}{\log n}\right)}\right) \\ &\quad \times \left(1 + O\left(\frac{\log\log n}{\sqrt{\log n}}\right)\right). \end{aligned}$$

Note that by Taylor's expansion, the right-hand side of (2.46) is of the form

$$\log j\left(1 + O\left(\frac{\log\log n}{\log k_n}\right)\right)\left(1 + O\left(\frac{\log\log n}{\sqrt{\log n}}\right)\right).$$

Hence, if we define $\alpha_{n,j}$ by

$$\alpha_{n,j}\log j := (1 - \varepsilon_n)2\log n\left(1 - \sqrt{1 - \frac{\log j}{\log n}}\right),$$

with $\varepsilon_n$ well chosen, then we get

$$\mathbb{P}\left[\left(\bigcup_{k_n \leq j \leq n\tilde{\lambda}_n} \sqrt{2\log n}\,(X_{j,n} - X_{1,n}) \leq \alpha_{n,j}\log j\right) \cap \Omega_{k_n}\right] = \mathbb{P}(\varnothing) = 0$$

and

$$(2.47) \qquad \begin{aligned} &\mathbb{P}\left[\left(\bigcup_{k_n \leq j \leq n\tilde{\lambda}_n} \sqrt{2\log n}\,(X_{j,n} - X_{1,n}) \leq \alpha_{n,j}\log j\right) \cap \Omega_{k_n}^c\right] \\ &\qquad \leq \mathbb{P}[\Omega_{k_n}^c] \leq 4\left(\frac{1}{(\log n)^{1+\delta}} + \frac{e^{-k_n/16}}{1 - e^{-1/16}}\right), \end{aligned}$$

and therefore (2.21).

To prove (2.25), that is, when $j \in J_{k_n} := \{j \in [k_n, n\tilde{\lambda}_n), \text{ such that } \log j/\log n \downarrow 0\}$, by Taylor's expansion in (2.41), we have on $\Omega_{n,\delta}$,

$$(2.48) \qquad \begin{aligned} &\sqrt{2\log n}\,(X_{j,n} - X_{1,n}) \\ &= \left(-\log\left(\frac{U_{1,n}}{U_{j,n}}\right) + o\left(\log\left(\frac{U_{1,n}}{U_{j,n}}\right)\right)\right)\left(1 + O\left(\frac{\log\log n}{\sqrt{\log n}}\right)\right). \end{aligned}$$

Therefore, to evaluate

$$\mathbb{P}\left[\left(\bigcup_{j \in J_{k_n}} \sqrt{2\log n}(X_{j,n} - X_{1,n}) \leq \alpha^*_{j,\eta} \log j\right) \cap \Omega_{n,\delta}\right],$$

where $\alpha^*_{j,\eta} := 1 - 1/(\log j)^{1-\eta}$, it is enough to estimate

$$\mathbb{P}_{k_n} := \mathbb{P}\left[\left(\bigcup_{j \in J_{k_n}} -\log\left(\frac{U_{1,n}}{U_{j,n}}\right) \leq \alpha^*_{j,\eta} \log j\right) \cap \Omega_{n,\delta}\right].$$

By using (2.13) with $c$ such that $(1-\varepsilon)c/4 > 1$ and $j_0 = k_n$, and by defining

$$\Omega^*_{k_n} := \bigcap_{j \geq k_n} \left\{\log\left(\frac{U_{1,n}}{U_{j,n}}\right) \leq -Z + \sqrt{\frac{c\log j}{j}} - \log j\right\},$$

we get that $\mathbb{P}[\Omega^*_{k_n}] \geq 1 - 1/(k_n^{(1-\varepsilon)c/4-1})$. Now it is easy to check that

$$\mathbb{P}_{k_n} \leq \mathbb{P}\left[\bigcup_{j \in J_{k_n}} \left\{Z \leq \left(\sqrt{\frac{c}{j\log j}} + (\alpha^*_{j,\eta} - 1)\right)\log j\right\}\right] + \mathbb{P}[(\Omega^*_{k_n})^c]$$

$$\leq \mathbb{P}[Z \leq -|x_{k_n}|\log k_n] + \mathbb{P}[(\Omega^*_{k_n})^c]$$

with $|x_j| := 1/(\log j)^{1-\eta} - \sqrt{c/j\log j}$, where we have used at the last step that $(|x_j|\log j)_{j \geq k_n}$ is an increasing function of $j$. Hence via (2.11), we get

$$(2.49) \quad \mathbb{P}_{k_n} \leq \exp\left\{-\frac{\sqrt{3}}{2}(\log k_n)^\eta + \frac{\sqrt{3}}{2}\sqrt{\frac{c\log k_n}{k_n}} + \frac{15}{4}\right\} + \frac{1}{k_n^{(1-\varepsilon)c/4-1}},$$

which implies (2.25). $\quad\square$

PROPOSITION 2.3. *For all $\delta > 0$, for all $k_n$ satisfying (2.17) and for $\Omega^*_n \subset \Omega$ given by $\Omega^*_n := \Omega_{n,k} \cap \Omega_{n,\delta} \cap \tilde{\Omega}_{k_n}$ defined, respectively, in (2.52), (2.40) and (2.43), we have*

$$(2.50) \qquad \mathbb{P}[\Omega^*_n] \geq 1 - \frac{6}{(\log n)^{1+\delta}},$$

*and on $\Omega^*_n$, we have, $\forall j, \ 1 \leq j \leq k_n$,*

$$(2.51) \qquad \sqrt{2\log n}(X_{j,n} - X_{1,n}) = \sum_{i=1}^{j} \frac{W_i}{i} + O\left(\frac{\log(k_n)\log\log n}{\log n}\right).$$

PROOF. For $\gamma > 0$ let us define

$$(2.52) \qquad \Omega_{n,k} := \left\{\sup_{1 \leq j \leq k_n} U_{j,n} = U_{k_n,n} \leq (1+\gamma)\frac{k_n}{n}\right\}.$$

It follows from (2.3) and (2.9) that $\mathbb{P}[\Omega_{n,k}^c] \le 4\exp(-k_n(\gamma^2/(4(2+\gamma)^2)))$.

From now on, our work space will be $\hat{\Omega}_{n,k}$. By using (2.1), we can write

$$\sqrt{2\log n}(X_{j,n} - X_{1,n})$$
$$= \sqrt{2\log n}\sqrt{2\log\left(\frac{1}{U_{1,n}}\right)}$$
$$\times \Bigg\{ 1 - \sqrt{1 + \frac{\log(U_{1,n}/U_{j,n})}{\log(1/U_{1,n})}}$$
$$- \frac{\log\log(1/U_{1,n})}{4\log(1/U_{1,n})} + \frac{\log\log(1/U_{j,n})}{4\sqrt{\log(1/U_{j,n})}\sqrt{\log(1/U_{1,n})}}$$
$$- \frac{\log(4\pi)}{4}\left(\frac{1}{\log(1/U_{1,n})} - \frac{1}{\sqrt{\log(1/U_{j,n})}\sqrt{\log(1/U_{1,n})}}\right)$$
$$+ O\left(\left(\frac{\log\log(1/U_{1,n})}{\log(1/U_{1,n})}\right)^2 + \frac{(\log\log(1/U_{j,n}))^2}{\log^{1/2}(1/U_{1,n})\log^{3/2}(1/U_{j,n})}\right)\Bigg\},$$

but by using the same type of arguments as in the proof of Proposition 2.2, Case (ii), we have, on $\Omega_{n,\delta} \cap \tilde{\Omega}_{k_n}$ [defined in (2.40) and (2.43), resp.],

$$\sup_{1 \le j \le k_n}\left(\frac{1}{\log(1/U_{1,n})} - \frac{1}{\sqrt{\log(1/U_{j,n})}\sqrt{\log(1/U_{1,n})}}\right)$$
$$= \frac{1}{\log(1/U_{1,n})}\left(1 - \frac{1}{\sqrt{(\log(1/U_{k_n,n}))/(\log(1/U_{1,n}))}}\right)$$
$$= O\left(\frac{\log k_n}{\log^2 n}\right),$$

$$\sup_{1 \le j \le k_n}\left(\frac{\log\log(1/U_{1,n})}{\log(1/U_{1,n})} - \frac{\log\log(1/U_{1,n})}{\sqrt{\log(1/U_{j,n})}\sqrt{\log(1/U_{1,n})}}\right)$$
$$= \frac{\log\log(1/U_{1,n})}{\log(1/U_{1,n})}\left(1 - \frac{1}{\sqrt{(\log(1/U_{k_n,n}))/(\log(1/U_{1,n}))}}\right)$$
$$= O\left(\frac{\log(k_n)\log\log n}{\log^2 n}\right)$$

and also

$$\sup_{1 \le j \le k_n}\left(\left(\frac{\log\log(1/U_{1,n})}{\log(1/U_{1,n})}\right)^2 + \frac{(\log\log(1/U_{j,n}))^2}{\log^{1/2}(1/U_{1,n})\log^{3/2}(1/U_{j,n})}\right)$$

$$= o\left(\frac{\log(k_n)\log\log n}{\log^2 n}\right).$$

Therefore, we obtain that, on $\Omega_{n,\delta}\cap\bar\Omega_{k_n}$, $\forall j, 1\le j\le k_n$,

$$
\begin{aligned}
&\sqrt{2\log n}(X_{j,n}-X_{1,n})\\
&\quad=2\sqrt{\log n\log\frac{1}{U_{1,n}}}\left\{1-\sqrt{1+\frac{\log(U_{1,n}/U_{j,n})}{\log(1/U_{1,n})}+O\left(\frac{\log(k_n)\log\log n}{\log^2 n}\right)}\right\}\\
&\quad=\sqrt{\log n\log\frac{1}{U_{1,n}}}\left(-\frac{\log(U_{1,n}/U_{j,n})}{\log(1/U_{1,n})}+O\left(\frac{\log(k_n)\log\log n}{\log^2 n}\right)\right)\\
&\quad=\sum_{i=1}^{j}\frac{W_i}{i}+O\left(\frac{\log(k_n)\log\log n}{\log n}\right)
\end{aligned}
$$

by using Taylor's expansion in the second equality, and (2.10) and (2.40) in the last equality. $\square$

**3. Applications to REM.** We are now able to answer the questions raised in the Introduction by using the tools developed in the previous section. We choose $n:=2^N$ and use $N$ instead of $n(N)$, since $N$ is related to the volume of our system. Otherwise, we keep the same notation used up to now, in particular $\beta_c=\sqrt{2\log 2}$, with $\beta/\beta_c>1$.

The Gibbs measure, defined in (1.4), can also be expressed as

$$(3.1)\qquad \mu_{N,\beta}(\Psi)=\sum_{i=1}^{2^N}\Psi((\sigma)_i)\mu_{N,\beta}((\sigma)_i)=\frac{\sum_{i=1}^{2^N}\Psi((\sigma)_i)e^{-\beta\sqrt{N}X_i}}{\sum_{i=1}^{2^N}e^{-\beta\sqrt{N}X_i}},$$

where $(X_i, i\in\{1,\dots,2^N\})$ are i.i.d. standard Gaussian random variables [associated to $\sigma=\{(\sigma)_i, i\in\{1,\dots,2^N\}\}$, an enumeration of $\{-1,+1\}^N$] and $\Psi$ is a real function defined on $\{-1,+1\}^N$, such that $\|\Psi\|_\infty:=\sup_{\sigma\in\{-1,+1\}^N}|\Psi(\sigma)|<\infty$.

We denote by $(\tilde\sigma)_j$ the $(\sigma)_i$ to which $X_i=X_{j,2^N}$ is associated, so we can write

$$(3.2)\qquad \mu_{N,\beta}(\Psi)=\frac{\sum_{i=1}^{2^N}\Psi((\tilde\sigma)_i)\exp(-(\beta/\beta_c)\sqrt{2\log 2^N}X_{i,2^N})}{\sum_{i=1}^{2^N}\exp(-(\beta/\beta_c)\sqrt{2\log 2^N}X_{i,2^N})}.$$

3.1. *Approximation of the Gibbs measure.* The first result of this section is related to how many terms of the order statistics we need to ensure a good approximation of the Gibbs measure.

THEOREM 3.1.  *Let $k_N$ satisfy*

$$(3.3) \qquad k_N \uparrow \infty, \qquad \frac{k_N}{N} \uparrow \infty \quad and \quad \frac{\log k_N}{N} \downarrow 0 \qquad as \ N \to \infty.$$

*There exist $\Omega_N \subset \Omega$ and $N_0$ such that for all $N > N_0$,*

$$(3.4) \qquad \mathbb{P}[\Omega_N] \geq 1 - \frac{4}{(N \log 2)^{1+\delta}}$$

*and on $\Omega_N$ we have*

$$(3.5) \qquad \begin{aligned} &\mu_{N,\beta}(\Psi) \\ &= \frac{\sum_{i=1}^{k_N} \Psi((\tilde{\sigma})_i) \exp(-(\beta/\beta_c)\sqrt{2 \log 2^N}(X_{i,2^N} - X_{1,2^N})) + B_N(\Psi)}{\sum_{i=1}^{k_N} \exp(-(\beta/\beta_c)\sqrt{2 \log 2^N}(X_{i,2^N} - X_{1,2^N})) + B_N(1)} \end{aligned}$$

*with*

$$(3.6) \qquad |B_N(\Psi)| \leq \frac{2\|\Psi\|_\infty}{\beta/\beta_c - 1} \frac{1}{(k_N - 1)^{(\beta/\beta_c - 1)/2}}.$$

PROOF.  From (3.2), we deduce (3.5) with $B_N(\Psi) := \sum_{i=k_N+1}^{2^N} \Psi((\tilde{\sigma})_i) \times \exp(-(\beta/\beta_c)\sqrt{2 \log 2^N}(X_{i,2^N} - X_{1,2^N}))$. We have to prove only that $B_N$ satisfies (3.6), which is a consequence of Proposition 2.2.

Let $\varphi_N(i) := \exp(-(\beta/\beta_c)\sqrt{2 \log 2^N}(X_{i,2^N} - X_{1,2^N}))$ and let $\tilde{\lambda}_N := \tilde{\lambda}_{2^N}$ satisfy (2.19). Then

$$B_N(\Psi) \leq \|\Psi\|_\infty \left( \sum_{i=k_N}^{2^N \tilde{\lambda}_N} \varphi_N(i) + \sum_{i=2^N \tilde{\lambda}_N + 1}^{2^N} \varphi_N(i) \right).$$

On one hand, by using (2.24), we have on $\Omega_N$,

$$(3.7) \qquad \begin{aligned} \sum_{i=2^N \tilde{\lambda}_N + 1}^{2^N} \varphi_N(i) &\leq \exp\left\{ -\left( \frac{\beta}{\beta_c}(1-\varepsilon) - 1 \right) \log 2^N \right\} = 2^{-N\tau} \\ &\text{with } \tau := \frac{\beta}{\beta_c}(1-\varepsilon) - 1 > 0. \end{aligned}$$

On the other hand, by using (2.21) with $\varepsilon = (1 - \beta_c/\beta)/2$, we have, on $\Omega_n$,

$$\sum_{k_N \leq i \leq 2^N \tilde{\lambda}_N} \varphi_N(i)$$

(3.8)
$$\leq \sum_{k_N \leq i \leq 2^N \tilde{\lambda}_N} \exp\left\{-\frac{\beta}{\beta_c}(1-\varepsilon)\frac{2\log i}{1+\sqrt{1-\log i/\log 2^N}}\right\}$$

$$\leq \sum_{k_N \leq i \leq 2^N \tilde{\lambda}_N} \frac{1}{i^{(1+\tau)(1-\varepsilon)}}$$

$$\leq \frac{2}{\beta/\beta_c - 1} \frac{1}{(k_N-1)^{(\beta/\beta_c-1)/2}}.$$

Hence (3.7) and (3.8) imply (3.6). $\square$

REMARK. An example of $k_N$ that satisfies (3.3) is $k_N = N\log_p(N)$, where $\log_p = \log\log_{p-1}$.

Now we can define the random probability measure on $\{-1, +1\}^N$ by

(3.9) $\mu_{k_N,\beta}^{(1)}(\Psi) := \dfrac{\sum_{k=1}^{k_N} \Psi((\tilde{\sigma})_k) \exp(-(\beta/\beta_c)\sqrt{2\log 2^N}(X_{i,2^N} - X_{1,2^N}))}{\sum_{k=1}^{k_N} \exp(-(\beta/\beta_c)\sqrt{2\log 2^N}(X_{i,2^N} - X_{1,2^N}))}$,

and thus show that this measure is a good approximation of the Gibbs measure. Indeed, by using the total variation distance between two measures given by

(3.10) $d_{\mathrm{TV}}(\mu, \nu) = \sup_{\{\Psi : \|\Psi\|_\infty = 1\}} |\mu(\Psi) - \nu(\Psi)|$,

the first Borel–Cantelli lemma and the fact that the denominator in (3.5) can be bounded from below by 1, it is immediate to check the following corollary:

COROLLARY 3.1. For all $\beta$ such that $\beta/\beta_c > 1$, for all $k_N$ that satisfy (3.3), we have with probability 1, for all but a finite number of indices $N$,

(3.11) $d_{\mathrm{TV}}(\mu_{N,\beta}, \mu_{k_N,\beta}^{(1)}) \leq \dfrac{4}{\beta/\beta_c - 1} \dfrac{1}{(k_N-1)^{(\beta/\beta_c-1)/2}}.$

With this last result and Proposition 2.3 in hand, we can define a second approximation of the Gibbs measure by the random measure on $\{-1, +1\}^N$,

(3.12) $\mu_{k_N,\beta}^{(2)}(\Psi) := \dfrac{\sum_{k=1}^{k_N} \Psi((\tilde{\sigma})_k) \exp(-\beta/\beta_c \sum_{\ell=1}^{k} W_\ell/\ell)}{\sum_{k=1}^{k_N} \exp(-\beta/\beta_c \sum_{\ell=1}^{k} W_\ell/\ell)}$,

where $(W_\ell, 1 < \ell < k_N)$ is a family of i.i.d. exponential random variables with mean 1.

PROPOSITION 3.1. *There exists an absolute constant $C$ such that, for all $\beta$ satisfying $\beta/\beta_c > 1$, for all $k_N$ satisfying*

$$(3.13) \quad k_N \uparrow \infty, \qquad \frac{k_N}{N} \uparrow \infty \quad and \quad \frac{(\log k_N)(\log N)}{N} \downarrow 0 \qquad as \ N \to \infty,$$

*we have with probability 1, for all but a finite number of indices $N$:*

$$
\begin{aligned}
& d_{\mathrm{TV}}(\mu_{k_N,\beta}, \mu^{(2)}_{k_N,\beta}) \\
(3.14) \quad & \leq C \frac{\beta}{\beta_c} \left( \frac{(\log k_N)(\log N)}{N} \right) \exp\left\{ C \frac{\beta}{\beta_c} \frac{(\log k_N)(\log N)}{N} \right\} \\
& \quad + \frac{4}{\beta/\beta_c - 1} \frac{1}{(k_N - 1)^{(\beta/\beta_c - 1)/2}}.
\end{aligned}
$$

PROOF. We have

$$
\begin{aligned}
& d_{\mathrm{TV}}(\mu^{(1)}_{k_N,\beta}, \mu^{(2)}_{k_N,\beta}) \\
(3.15) \quad & \leq C \frac{\beta}{\beta_c} \left( \frac{(\log k_N)(\log N)}{N} \right) \exp\left\{ C \frac{\beta}{\beta_c} \frac{(\log k_N)(\log N)}{N} \right\}
\end{aligned}
$$

as a consequence of Proposition 2.3, the fact that $\forall x, |e^x - 1| \leq |x| e^{|x|}$ and the first Borel–Cantelli lemma. Applying the triangle inequality, (3.11) and (3.15) gives (3.14). $\square$

REMARK. Note that (3.13) is satisfied if we choose $k_N = N \log_p(N)$, where $\log_p(x) = \log \log_{p-1}(x)$ and $\log_1(x) = \log x$.

3.2. *Some properties of the Gibbs measure.* The previous results help us to get more information on the support of the Gibbs measure. Indeed:

PROPERTY 1. *As an immediate consequence of Corollary 3.1 and Proposition 3.1, we have: for $k_N$ satisfying (3.13),*

$$(3.16) \quad \limsup_{N \uparrow \infty} \mu_{N,\beta}(\{(\sigma)_1, \ldots, (\sigma)_{k_N}\}^c) = 0 \qquad a.s.$$

PROPERTY 2. *When $\beta_c < \beta < \infty$, the Gibbs measure cannot be concentrated on the minimum, that is, $\mu_{k_N,\beta} \neq \delta_{(\tilde{\sigma})_1}$.*

PROOF. Considering (3.12) with the trial function $\Psi((\tilde{\sigma})_i) = \mathbb{1}_{\{i=1\}} \forall i \geq 1$, we can write

$$d_{\mathrm{TV}}(\mu^{(2)}_{k_N,\beta}, \delta_{(\tilde{\sigma})_1})$$

$$\geq \frac{\sum_{k=2}^{k_N} \exp(-\beta/\beta_c \sum_{\ell=1}^{k} W_\ell/\ell)}{1 + \sum_{k=2}^{k_N} \exp(-\beta/\beta_c \sum_{\ell=1}^{k} W_\ell/\ell)}$$

$$\geq \frac{\exp(-\beta/\beta_c(W_1 + W_2/2))}{1 + \exp(-\beta/\beta_c(W_1 + W_2/2))},$$

which in combination with Proposition 3.1, provides that

$$\liminf_{N\uparrow\infty} d_{\mathrm{TV}}(\mu_{N,\beta}, \delta_{(\tilde{\sigma})_1}) \geq \frac{\exp(-\beta/\beta_c(W_1 + W_2/2))}{1 + \exp(-\beta/\beta_c(W_1 + W_2/2))} > 0. \qquad \square$$

To go on with the properties, we need the following lemma.

LEMMA 3.1. *For all $\beta > 0$,*

$$(3.17) \qquad \lim_{k\to\infty} \sum_{j=2}^{k} \exp\left(-\beta/\beta_c \sum_{\ell=1}^{j} W_\ell/\ell\right) = \zeta(\beta) \qquad a.s.$$

*If $\beta > \beta_c$, then $\mathbb{E}[\zeta(\beta)] < \infty$. In particular, $\mathbb{P}[0 < \zeta(\beta) < \infty] = 1$. Moreover, for any given sequence $\beta_p$ such that $\lim_{p\uparrow\infty} \beta_p = \infty$, there exists a subset $\tilde{\Omega} \subset \Omega$ such that $\mathbb{P}[\tilde{\Omega}] = 1$ and, on $\tilde{\Omega}$,*

$$(3.18) \qquad \lim_{p\uparrow\infty} \zeta(\beta_p) = 0.$$

PROOF. Let $\Sigma_k$ be the $\sigma$-algebra generated by $(W_1, \ldots, W_k)$ and let

$$\zeta_k(\beta) = \sum_{j=2}^{k} \exp\left(-\frac{\beta}{\beta_c} \sum_{\ell=1}^{j} \frac{W_\ell}{\ell}\right).$$

The positive random variable $\zeta_k(\beta)$ is a supermartingale, since we have

$$(3.19) \qquad \mathbb{E}[\zeta_{k+1}(\beta)|\Sigma_k] = \frac{1}{1 + \beta/(\beta_c(k+1))} \zeta_k(\beta) \leq \zeta_k(\beta).$$

Therefore, the increasing sequence $\zeta_k(\beta)$ converges almost surely to some random variable $\zeta(\beta)$. On the other hand, using that $\forall x \geq 0$, $\log(1+x) \geq x - (x^2)/2$, $\sum_{\ell=1}^{j}(\ell)^{-1} \geq \log j$ and $\sum_{\ell=1}^{j}(\ell)^{-2} \leq 2 - 1/j \leq 2$ provides

$$\mathbb{E}[\zeta_k(\beta)] = \sum_{j=2}^{k} \prod_{\ell=1}^{j} \frac{1}{1 + \beta/(\beta_c k)}$$

$$(3.20) \qquad \leq e^{\beta^2/\beta_c^2} \sum_{j=2}^{k} \frac{1}{j^{\beta/\beta_c}}$$

$$\leq e^{\beta^2/\beta_c^2} \sum_{j=2}^{\infty} \frac{1}{j^{\beta/\beta_c}} < \infty,$$

where at the last step we have used $\beta/\beta_c > 1$. Hence, by the Fatou–Lebesgue lemma, we get $\mathbb{E}[\zeta(\beta)] < \infty$ and so $\mathbb{P}[0 < \zeta(\beta) < \infty] = 1$.

Proving (3.18) requires a little more care, since the probability subset where $\zeta_k(\beta)$ converges, denoted by $\Omega(\beta)$, depends a priori on $\beta$.

We can write $\zeta_k(\beta) = \exp(-(\beta/\beta_c)W_1)\tilde{\zeta}_k(\beta)$ with

$$\tilde{\zeta}_k(\beta) = \sum_{j=2}^{k} \exp\left(-\frac{\beta}{\beta_c}\sum_{\ell=2}^{j}\frac{W_\ell}{\ell}\right).$$

By the same kind of computations we did previously, we can show that, for any sequence $\beta_c < \beta_p \uparrow \infty$, there exist some positive integrable random variables $\tilde{\zeta}(\beta_p)$ such that on the subspace $\tilde{\Omega} := \bigcap_{\beta_p}\Omega(\beta_p)$ of $\mathbb{P}$ probability 1, we have $\lim_{k\uparrow\infty}\tilde{\zeta}_k(\beta_p) = \tilde{\zeta}(\beta_p)$.

Since for all $k$, $\tilde{\zeta}_k(\beta_q) \leq \tilde{\zeta}_k(\beta_p)$ when $\beta_p \leq \beta_q$, we have, also on $\tilde{\Omega}$, $\tilde{\zeta}(\beta_q) \leq \tilde{\zeta}(\beta_p)$. Therefore, denoting $\beta^* = \inf\{\beta_p > 1 + \beta_c\}$, since $\mathbb{P}[\tilde{\zeta}(\beta^*) < \infty] = 1$ and $\mathbb{P}[\lim_{\beta_p \uparrow \infty}\exp(-(\beta/\beta_c)W_1) = 0] = 1$, with probability 1 we get that

$$0 \leq \limsup_{\beta_p\uparrow\infty}\zeta(\beta_p) \leq \limsup_{\beta_p\uparrow\infty} e^{(-\beta/\beta_c)W_1}\tilde{\zeta}(\beta^*) = 0. \qquad \square$$

PROPERTY 3. *At zero temperature, the Gibbs measure is a.s. the point mass at the minimum.*

PROOF. We have

$$d_{\mathrm{TV}}(\mu_{k_N,\beta}^{(2)}, \delta_{(\tilde{\sigma})_1}) \leq 2\frac{\sum_{k=2}^{k_N}\exp(-\beta/\beta_c\sum_{\ell=1}^{k}W_\ell/\ell)}{1 + \sum_{k=2}^{k_N}\exp(-\beta/\beta_c\sum_{\ell=1}^{k}W_\ell/\ell)},$$

which gives with probability 1, via Lemma 3.1,

$$\limsup_{N\uparrow\infty}d_{\mathrm{TV}}(\mu_{k_N,\beta}^{(2)}, \delta_{(\tilde{\sigma})_1}) \leq \frac{2\zeta(\beta)}{1 + \zeta(\beta)}.$$

Then Proposition 3.1 provides with probability 1 that

$$\limsup_{N\uparrow\infty}d_{\mathrm{TV}}(\mu_{N,\beta}, \delta_{(\tilde{\sigma})_1}) \leq \frac{2\zeta(\beta)}{1 + \zeta(\beta)}.$$

Now by using (3.18), we get with probability 1,

$$\limsup_{\beta_p\uparrow\infty}\limsup_{N\uparrow\infty}d_{\mathrm{TV}}(\mu_{N,\beta}, \delta_{(\tilde{\sigma})_1}) = 0. \qquad \square$$

REMARK. Note that there is a result in [29], in the case of "REM for size dependence," that deals with the quality of the concentration of the Gibbs measure on the minimizer of the energy when the temperature goes to zero with the system size.

PROPERTY 4. *With probability 1, the Gibbs measure is not the uniform measure on the $k_N$ first minima $\{(\tilde{\sigma})_1, \ldots, (\tilde{\sigma})_{k_N}\}$.*

PROOF. Let $\nu_{k_N}$ be the uniform measure on $k_N$ points. Considering (3.12) with $\Psi((\tilde{\sigma})_i) = \mathbb{1}_{(i=1)} \forall i \geq 1$, we can write

$$d_{\mathrm{TV}}(\mu_{k_N,\beta}^{(2)}, \nu_{k_N}) \geq \frac{1}{1 + \sum_{k=2}^{k_N} \exp(-\beta/\beta_c \sum_{\ell=1}^{k} W_\ell/\ell)} - \frac{1}{k_N},$$

which, combined with Lemma 3.1 and Proposition 3.1, gives that with probability 1,

$$\liminf_{N \uparrow \infty} d_{\mathrm{TV}}(\mu_{N,\beta}, \nu_{k_N}) \geq \frac{1}{1 + \zeta(\beta)} > 0. \qquad \square$$

PROPERTY 5. *With probability 1, the Gibbs measure is not the measure $\mu_{k_0,\beta}$ for a finite $k_0$.*

PROOF. Because of Proposition 3.1, it is enough to prove that for a finite $k_0$, $\liminf_{N \uparrow \infty} d_{\mathrm{TV}}(\mu_{k_0,\beta}^{(2)}, \mu_{k_N,\beta}^{(2)})$ is bounded from below by a quantity that does not go to zero. Considering (3.12) with $\Psi$ defined by $\Psi(\sigma \in \{(\tilde{\sigma})_1, \ldots, (\tilde{\sigma})_{k_0}\}) = 0$ and $\Psi(\sigma = (\tilde{\sigma})_i) = 1$ for all $k_0 \leq i \leq k_N$, we get

$$(3.21) \quad \begin{aligned} d_{\mathrm{TV}}(\mu_{k_0,\beta}^{(2)}, \mu_{k_N,\beta}^{(2)}) &\geq \mu_{k_N,\beta}^{(2)}(\{(\tilde{\sigma})_1, \ldots, (\tilde{\sigma})_{k_0}\}^c) \\ &\geq \frac{\exp(-\beta/\beta_c \sum_{\ell=1}^{k_0+1} W_\ell/\ell)}{1 + \sum_{k=2}^{k_N} \exp(-\beta/\beta_c \sum_{\ell=1}^{k} W_\ell/\ell)}, \end{aligned}$$

and so, via Lemma 3.1, with probability 1,

$$\liminf_{N \uparrow \infty} d_{\mathrm{TV}}(\mu_{k_0,\beta}^{(2)}, \mu_{k_N,\beta}^{(2)}) \geq \frac{\exp(-\beta/\beta_c \sum_{\ell=1}^{k_0} W_\ell/\ell)}{1 + \zeta(\beta)} > 0. \qquad \square$$

PROPERTY 6. *With probability 1, for $\beta_c < \beta < \infty$, the Gibbs measure does not have all its mass concentrated on $\{(\tilde{\sigma})_1, \ldots, (\tilde{\sigma})_{k_0}\}$ for a finite $k_0$.*

PROOF. This result comes from $\mu_{k_0,\beta}[\{(\tilde{\sigma})_1, \ldots, (\tilde{\sigma})_{k_0}\}] = 1$ and

$$\begin{aligned} \mu_{N,\beta}&[\{(\tilde{\sigma})_1, \ldots, (\tilde{\sigma})_{k_0}\}] \\ &= \mu_{N,\beta}[\{(\tilde{\sigma})_1, \ldots, (\tilde{\sigma})_{k_0}\}] - \mu_{k_N,\beta}[\{(\tilde{\sigma})_1, \ldots, (\tilde{\sigma})_{k_0}\}] \\ &\quad + \mu_{k_N,\beta}[\{(\tilde{\sigma})_1, \ldots, (\tilde{\sigma})_{k_0}\}] - \mu_{k_0,\beta}[\{(\tilde{\sigma})_1, \ldots, (\tilde{\sigma})_{k_0}\}] + 1 \end{aligned}$$

if one notes that

$$\begin{aligned} \mu_{k_N,\beta}&[\{(\tilde{\sigma})_1, \ldots, (\tilde{\sigma})_{k_0}\}] - \mu_{k_0,\beta}[\{(\tilde{\sigma})_1, \ldots, (\tilde{\sigma})_{k_0}\}] \\ &= \mu_{k_0,\beta}[\{(\tilde{\sigma})_1, \ldots, (\tilde{\sigma})_{k_0}\}^c] - \mu_{k_N,\beta}[\{(\tilde{\sigma})_1, \ldots, (\tilde{\sigma})_{k_0}\}^c] \\ &= -\mu_{k_N,\beta}[\{(\tilde{\sigma})_1, \ldots, (\tilde{\sigma})_{k_0}\}^c]. \end{aligned}$$

By using (3.21) and Proposition 3.1, we get

$$\limsup_{N\uparrow\infty} \mu_{N,\beta}[\{(\tilde{\sigma})_1,\ldots,(\tilde{\sigma})_{k_0}\}]$$

$$\leq 1 - \frac{\exp(-\beta/\beta_c \sum_{\ell=1}^{k_0+1} W_\ell/\ell)}{1+\zeta(\beta)} < 1. \qquad \square$$

3.3. *Representation of the Gibbs measure.* Thanks to the two approximations of the Gibbs measure given in (3.9) and (3.12), we are now able to propose a representation of the Gibbs measure and a way to simulate it. First we consider the measure $\mu_{k_N,\beta}^{(1)}$ defined in (3.9), that we can simulate in the following way. [Note that we have to pay attention to the relationship between the $(\tilde{\sigma}_i)$ and the $(\sigma_i)$, since the $(\tilde{\sigma}_i)$ are needed for the simulation.] Let $(U_1,\ldots,U_{2^N})$ denote a $2^N$ sample of uniformly distributed random variables of order statititics $(U_{1,2^N} \leq \cdots \leq U_{2^N,2^N})$. We want to construct $\exp\{-(\beta/\beta_c)\sqrt{2\log 2^N}(X_{k,2^N}-X_{1,2^N})\}$, because of (2.2). We consider and simulate only the $k_N$ last terms $U_{2^N,2^N},\ldots,U_{2^N-k_N+1,2^N}$, with $k_N$ satisfying (3.3). (Note that this can be done without having first to order uniform variates by using some faster algorithm (cf., e.g., [23], Section 3.26).) Then *independently*, we choose, one after the other, $k_N$ spin configurations in $\{-1,+1\}^N$ without replacement. This defines an ordered sequence of configurations that we call $((\sigma)_1,\ldots,(\sigma)_{k_N})$. We make the following claim.

CLAIM 1. *The random measure defined by*

$$\tilde{\mu}_{k_N,\beta}^{(1)}(\Psi)$$

(3.22)
$$:= \left\{ \Psi((\sigma)_1) \right.$$
$$\left. + \sum_{k=2}^{k_N} \Psi((\sigma)_k) \exp\left(-(\beta/\beta_c)\sqrt{2\log 2^N}(G(U_{1,2^N})-G(U_{k,2^N}))\right) \right\}$$
$$\times \left\{ 1 + \sum_{k=2}^{k_N} \exp\left(-(\beta/\beta_c)\sqrt{2\log 2^N}(G(U_{1,2^N})-G(U_{k,2^N}))\right) \right\}^{-1}$$

*has the same distribution as $\mu_{k_N,\beta}^{(1)}$.*

Indeed, the claim is an immediate consequence of the following technical lemma on spacings of $n$-independent uniform random variables. For $1 \leq j \leq n$, let $\ell_j$ be the index of the (almost surely unique) $U_i$ such that $U_{\ell_j} = U_{j,n}$.

LEMMA 3.2.   *For all $(x_i)_{1 \leq i \leq k}$ and for all integer $k$, we have*

$$\mathbb{P}[U_{i,n} \leq x_i, U_{i,n} = U_{l_i}, 1 \leq i \leq k]$$
$$= \mathbb{P}[U_{i,n} \leq x_i, 1 \leq i \leq k] \times \mathbb{P}[U_{i,n} = U_{l_i}, 1 \leq i \leq k].$$

*Moreover,*

$$\mathbb{P}[U_{i,n} = U_{l_i}, 1 \leq i \leq k] = \frac{(n-k)!}{n!}.$$

PROOF.   We have (cf., e.g., [13])

$$\mathbb{P}[U_{i,n} \leq x_i, 1 \leq i \leq k]$$
$$= \frac{n!}{(n-k)!} \int_0^{x_1} \cdots \int_0^{x_k} (1 - y_k)^{n-k} \mathbb{1}_{(0 \leq y_1 < y_2 < \cdots < y_k \leq 1)} \prod_{i=1}^k dy_i.$$

Moreover, we can show by induction that

$$\mathbb{P}[U_{i,n} \leq x_i, U_{i,n} = U_{l_i}, 1 \leq i \leq k]$$
$$= \int_0^{x_1} \cdots \int_0^{x_k} (1 - y_k)^{n-k} \mathbb{1}_{(0 \leq y_1 < y_2 < \cdots < y_k \leq 1)} \prod_{i=1}^k dy_i,$$

which completes the proof.   □

Thus we have the first representation of the Gibbs measure that can be simulated. An alternative and easy way to proceed is to consider $\mu_{k_N,\beta}^{(2)}$ defined in (3.12) instead of $\mu_{k_N,\beta}^{(1)}$. To simulate $\mu_{k_N,\beta}^{(2)}$, we just have to consider a $k_N$ sample of independent uniformly distributed random variables $U_1, \ldots, U_{k_N}$, with $k_n$ satisfying (3.13), then choose as before the $k_N$ spin configurations, $(\sigma)_1, \ldots, (\sigma)_{k_N}$ to construct the resulting measure

$$(3.23) \quad \tilde{\mu}_{k_N,\beta}^{(2)}(\Psi) := \frac{\Psi((\sigma)_1) + \sum_{k=2}^{k_N} \Psi((\sigma)_k) \exp(+\beta/\beta_c \sum_{l=1}^k (\log U_\ell)/\ell)}{1 + \sum_{k=2}^{k_N} \exp(+\beta/\beta_c \sum_{l=1}^k (\log U_\ell)/\ell)}.$$

Then we have, in the same way as for the above claim, that $\tilde{\mu}_{k_N,\beta}^{(2)} =^d \mu_{k_N,\beta}^{(2)}$.

Thus this second procedure needs only two independent samples: one of $k_N$ spin configurations chosen without replacement within $2^N$ and one of $k_N$ independent uniform random variables on $(0,1)$, so we do not have to work with the $k_n$ largest order statistics.

Note that, as mentioned previously, $k_N$ must satisfy (3.13) and so could be chosen as $N \log \log \log N$, for instance.

**Acknowledgments.**  We wish to thank the referees for their comments and for providing some recent references.

U.F.R. DE MATHÉMATIQUES
   ET INFORMATIQUE
UNIVERSITÉ RENÉ DESCARTES, PARIS V
45 RUE DES SAINTS-PÈRES
75270 PARIS CEDEX 06
FRANCE
E-MAIL: Marie.Kratz@math-info.univ-paris5.fr

CENTRE DE PHYSIQUE THÉORIQUE
CNRS LUMINY
CASE 907
13288 MARSEILLE CEDEX 09
FRANCE
E-MAIL: Pierre.Picco@cpt.univ-mrs.fr